\documentclass[MSNbibl,number,citesort,seceqn,dvips]{arxbj}
\usepackage{upgreek}
\usepackage{graphicx}

\aid{0}
\volume{20}
\issue{1}
\pubyear{2014}
\firstpage{231}
\lastpage{244}
\doi{10.3150/12-BEJ483} %kopijuoti is PTS

\makeatletter
\renewcommand{\pi}{\uppi}
\def\cal{\mathcal}
\def\RR{\mathbb{R}}
\def\argmax{\mathop{\operatorname{argmax}}}

\newtheorem{theorem}{Theorem}[section]
\newtheorem{corollary}{Corollary}[section]
\newtheorem{proposition}{Proposition}[section]
\newtheorem{conjecture}{Conjecture}[section]
\makeatother

\begin{document}
\begin{frontmatter}

\title{Chernoff's density is log-concave}
\runtitle{Chernoff's density is log-concave}

\begin{aug}
%%%% inicialai - be tarpu
\author[1]{\fnms{Fadoua} \snm{Balabdaoui}\thanksref{1}\ead[label=e1]{fadoua@ceremade.dauphine.fr}}%
\and
\author[2]{\fnms{Jon A.} \snm{Wellner}\corref{}\thanksref{2}\ead[label=e2]{jaw@stat.washington.edu}}
\runauthor{F. Balabdaoui and J.A. Wellner} %% auto
\address[1]{Centre de Recherche en Math\'ematiques de la D\'ecision,
Universit\'e Paris-Dauphine,
Paris, France. \printead{e1}}
\address[2]{Department of Statistics, University of Washington, Seattle,
WA 98195-4322, USA.\\ \printead{e2}}
\end{aug}

% HISTORY:
\received{\smonth{3} \syear{2012}}
\revised{\smonth{10} \syear{2012}}

% ABSTRACT
%
\begin{abstract}
We show that the density of $Z = \argmax\{ W(t) - t^2 \}$,
sometimes known as Chernoff's density, is log-concave.
We conjecture that Chernoff's density is strongly log-concave
or ``super-Gaussian'', and provide evidence in support of the conjecture.
\end{abstract}

% KEYWORDS
% visi is mazosios raides ir pagal abecele
%
\begin{keyword}
\kwd{airy function}
\kwd{Brownian motion}
\kwd{correlation inequalities}
\kwd{hyperbolically monotone}
\kwd{log-concave}
\kwd{monotone function estimation}
\kwd{Polya frequency function}
\kwd{Prekopa--Leindler theorem}
\kwd{Schoenberg's theorem}
\kwd{slope process}
\kwd{strongly log-concave}
\end{keyword}

\end{frontmatter}

%s1 #&#
\section{Introduction: Two limit theorems}
\label{sec:TwoLimitTheorems}

We begin by comparing two limit theorems.

First the usual central limit theorem:
Suppose that $X_1, \ldots, X_n $ are i.i.d. $EX_1 = \mu$, $E(X^2)<
\infty$, $\sigma^2 = \operatorname{Var}(X)$.
Then, the classical Central Limit theorem says that
\[
\sqrt{n} (\overline{X}_n - \mu) \rightarrow_d N \bigl(0,
\sigma^2 \bigr) .
\]
The Gaussian limit has density
\begin{eqnarray*}
 \phi_{\sigma} (x) &= &\frac{1}{\sqrt{2 \pi} \sigma} \exp \biggl(- \frac{x^2}{2 \sigma^2 }
\biggr) = \mathrm{e}^{-V(x)} ,
\\
 V(x) &=&- \log\phi_{\sigma} (x) = \frac{x^2}{2 \sigma^2} + \log ( \sqrt{2
\pi} \sigma),
\\
 V'' (x)& =& (- \log\phi_{\sigma}
)'' (x) = \frac{1}{\sigma^2} > 0 .
\end{eqnarray*}
Thus, $\log\phi_\sigma$ is concave, and hence $\phi_{\sigma}$ is a
\textit{log-concave density}.
As is well known, the normal distribution arises as a natural limit in
a wide range of
settings connected with sums of independent and weakly dependent random
variables;
see, for example, Le~Cam \cite{MR833276} %Le Cam
and Dehling and Philipp \cite{MR1958777}.

Now for a much less well-known limit theorem in the setting of monotone
regression.
Suppose that the real-valued function $r(x)$
is monotone increasing for $x \in[0,1]$.
For $i \in\{ 1, \ldots, n \}$, suppose that $x_i = i/(n+1)$,
$\varepsilon_i $ are i.i.d. with $E(\varepsilon_i)=0$,
$\sigma^2=E(\varepsilon_i^2 )< \infty$, and suppose that we observe
$(x_i , Y_i )$, $i=1,\ldots, n$, where
\[
Y_i = r(x_i) + \varepsilon_i \equiv
\mu_i + \varepsilon_i,\qquad  i \in\{ 1, \ldots, n \}.
\]
The isotonic estimator $\widehat{\underline{\mu}} $ of $\underline
{\mu} = (\mu_1 ,\ldots, \mu_n)$
is given by
\begin{eqnarray*}
 \widehat{\mu}_j& =& \max_{i \le j} \min_{k \ge j} \biggl
\{ \frac
{\sum_{l=i}^k Y_l}{k-i+1} \biggr\},
\\
\widehat{\underline{\mu}}& = &( \hat{\mu}_1, \ldots, \hat{\mu
}_n ) \equiv T \underline{Y}
\\
 &=& \mbox{least squares projection of } \underline{Y} \mbox { onto }
K_n,
\\
 K_n& = &\bigl\{ y \in\RR^n \dvt y_1 \le
\cdots\le y_n \bigr\} .
\end{eqnarray*}
For fixed $x_0 \in(0,1)$ with $x_j \le x_0 < x_{j+1}$ we set
$\hat{r}_n (x_0 ) \equiv\hat{r}_n (x_j ) = \hat{\mu}_j$.

Brunk \cite{MR0277070} %Brunk (1970)
showed that if
$r'(x_0)>0$ and if $r'$ is continuous in a neighborhood of $x_0$, then
\[
n^{1/3} \bigl( \widehat{r}_n (x_0) -
r(x_0) \bigr) \rightarrow_d \bigl(\sigma^2
r'(x_0)/2 \bigr)^{1/3} (2Z_1),
\]
where, with $\{ W (t) \dvt t \in\RR\}$ denoting a two-sided standard
Brownian motion
process started at $0$,
%
%e1.1 #&#
\begin{eqnarray}\label{BM-DescriptionOfZ2}
2Z_1 & = & \mbox{slope at zero of the greatest convex minorant of }
W(t) + t^2
\nonumber\\
& \stackrel{d} {=} & \mbox{slope at zero of the least concave majorant of } W(t)
- t^2
\\
& \stackrel{d} {=} & 2 \operatorname{argmin} \bigl\{ W(t) + t^2
\bigr\} .
\nonumber
\end{eqnarray}
The density $f$ of $Z_1$ is called {Chernoff's density}.
Chernoff's density appears in a number of nonparametric problems involving
estimation of a monotone function:
\begin{itemize}
\item Estimation of a monotone regression function $r$: see, for example,
Ayer \textit{et~al.} \cite{MR0073895}, %Ayer, Brunk, Ewing, Reid, Silverman (1955),
van Eeden~\cite{MR0083869}, %van Eeden (1957),
Brunk \cite{MR0277070}, %Brunk (1970),
and Leurgans \cite{MR642740}. %Leurgans (yy).
\item Estimation of a monotone decreasing density: see
Grenander \cite{MR0086459}, %Grenander (1956),
Prakasa~Rao \cite{MR0267677}, %Prakasa Rao (1969),
and Groeneboom \cite{MR822052}. %Groeneboom (1985).
\item Estimation of a monotone hazard function:
Grenander \cite{MR0093415}, %Grenander (1956),
Prakasa~Rao \cite{MR0260133}, %Prakasa Rao (1970),
Huang and Zhang \cite{MR1311975}, %Huang and Zhang (xx),
Huang and Wellner \cite{MR1334065}. % Huang and Wellner (yy),
\item Estimation of a distribution function with interval censoring:
Groeneboom and Wellner \cite{MR1180321}, %Groeneboom-Wellner (1992),
Groeneboom \cite{MR1600884}.
\end{itemize}
In each case:
\begin{itemize}
\item There is a monotone function $m$ to be estimated.
\item There is a natural nonparametric estimator $\widehat{m}_n $.
\item If $m'(x_0) \not= 0$ and $m'$ continuous at $x_0$, then
\[
n^{1/3} \bigl( \widehat{m}_n (x_0) -
m(x_0) \bigr) \rightarrow_d C (m, x_0)
2Z_1,
\]
where $2Z_1$ is as in (\ref{BM-DescriptionOfZ2}).
\end{itemize}
See Kim and Pollard \cite{MR1041391} %Kim and Pollard (1990)
for a unified approach to these types of problems.

The first appearance of $Z_1$ was in
Chernoff \cite{MR0172382}. %Chernoff (1964).
Chernoff \cite{MR0172382} considered estimation of the mode
of a (symmetric unimodal) density $f$ via the following simple estimator:
if $X_1,\ldots, X_n$ are i.i.d. with density $h$ and distribution
function $H$,
then for each fixed $a>0$ let
\[
\hat{x}_a \equiv \mbox{center of the interval of length } 2a
\mbox{ containing the most observations}.
\]
Let $x_a $ be the center of the interval of length $2a$ maximizing
$H(x+a) - H(x-a) = P( X \in(x-a,x+a])$. (Note that this $x_a$ is \textit{not} the mode if $f$ is
not symmetric.)
Then Chernoff shows:
\[
n^{1/3} (\hat{x}_a - x_a)
\rightarrow_d \biggl( \frac{ h(x_a +a)}{c} \biggr)^{1/3}
2Z_1,
\]
where
$c \equiv h' (x_a-a) - h'(x_a +a)$.
Chernoff also showed that the density $f_{Z_1}= f$ of $Z_1$ has the form
%
%e1.2 #&#
\begin{equation}
f(z) \equiv f_{Z_1}(z) = \tfrac{1}{2} g(z) g(-z),
\label{SymmetricProductFormChernoff}
\end{equation}
where
\[
g (t) \equiv\lim_{x \nearrow t^2} \frac{\partial}{\partial x} u (t,x) ,
\]
where, with $W$ standard Brownian motion,
\[
u(t,x) \equiv P^{(t,x)} \bigl(W(z) > z^2, \mbox{for some } z
\ge t \bigr)
\]
is a solution to the backward heat equation
\[
\frac{\partial}{\partial t} u(t,x) = - \frac{1}{2} \frac{\partial^2}{\partial x^2} u(t,x)
\]
under the boundary conditions
\[
u \bigl(t,t^2 \bigr) = \lim_{x \nearrow t^2} u(t,x) =1,\qquad
\lim_{x\rightarrow
-\infty} u(t,x) = 0 .
\]

Again let $W(t)$ be standard two-sided Brownian motion starting from
zero, and let $c>0$.
We now define
%
%e1.3 #&#
\begin{eqnarray}
Z_c \equiv\sup \bigl\{ t \in\RR\dvt W(t) - c t^2
\mbox{ is maximal} \bigr\} . \label{ArgmaxFunctional}
\end{eqnarray}
As noted above, $Z_c$ with $c=1$
arises naturally in the limit theory for nonparametric estimation of
monotone (decreasing) functions.
Groeneboom \cite{MR981568} (see also Daniels and Skyrme \cite{MR778595}) showed that for all $c>0$
the random variable
$Z_c$ has density
\[
f_{Z_c} (t) = \tfrac{1}{2} g_c (t) g_c
(-t),
\]
where $g_c$ has Fourier transform given by
%
%e1.4 #&#
\begin{equation}
\hat{g}_c (\lambda) = \int_{-\infty}^{\infty}
\mathrm{e}^{\mathrm{i}\lambda s} g_c (s) \,\mathrm{d}s = \frac{2^{1/3} c^{-1/3}}{\operatorname{Ai} (\mathrm{i} (2c^2)^{-1/3} \lambda)} .
\label{GroeneboomFT}
\end{equation}
Groeneboom and Wellner \cite{MR1939706} gave numerical computations of the density $f_{Z_1}$,
distribution function, quantiles, and moments.

Recent work on the distribution of the supremum
$M_c \equiv\sup_{t \in\RR} (W(t) - ct^2)$ is given in
Janson, Louchard and
  Martin-L{\"o}f \cite{JansonLouchard:10} and Groeneboom \cite{Groeneboom:10}.
Groeneboom \cite{Groeneboom:11} studies the number of vertices of the
greatest convex minorant
of $W(t) + t^2$ in intervals $[a,b]$ with $b-a\rightarrow\infty$; the
function $g_c$ with $c=1$ also plays
a key role there.

Our goal in this paper is to show that the density $f_{Z_c}$ is log-concave.
We also present evidence in support of the conjecture that $f_{Z_c}$ is
strongly log-concave:
that is, $(-\log f_{Z_c} )'' (t) \ge\mbox{some }  c>0$ for all $t \in
\RR$.

The organization of the rest of the paper is as follows:
log-concavity of $f_{Z_c}$ is proved in Section~\ref{sec:LogConcave}
where we also give graphical support for this property and present
several corollaries and related results.
In Section~\ref{sec:StrongLogConcave}, we give some partial results
and further graphical evidence
for strong log-concavity of $f \equiv f_{Z_1}$: that is,\vspace*{-1pt}
\[
(-\log f)'' (t) \ge(-\log f)''
(0) = 3.4052\ldots=1/(0.541912\ldots )^2 \equiv1/
\sigma_0^2\vspace*{-1pt}
\]
for all $t \in\RR$.
As will be shown in Section~\ref{sec:StrongLogConcave}, this is
equivalent to
$f(t) = \rho(t) \phi_{\sigma_0} (t)$ with $\rho$ log-concave.
In Section~\ref{sec:OpenProbs}, we briefly discuss some of the
consequences and corollaries
of log-concavity and strong log-concavity, sketch connections to some
results of
Bondesson \cite{MR1224674,MR1481175}, and list a few of the many further problems.\vspace*{-2pt}

%s2 #&#
\section{Chernoff's density is log-concave}
\label{sec:LogConcave}

Recall that a function $h$ is a \textit{P\'olya frequency function of
order $m$}
(and we write $h \in \mathrm{PF}_{m} $) if $K(x,y) \equiv h(x-y)$
is totally positive of order $m$: that is, $\operatorname{det} (H_m(\underline
{x}, \underline{y}) \ge 0$
for all choices of $x_1 \le\cdots\le x_m$ and $y_1 \le\cdots\le
y_m$ where
$H_m \equiv H_m (\underline{x}, \underline{y}) = (h(x_i- y_j))_{i,j=1}^m$.
It is well known and easily proved that a density $f$ is $\mathrm{PF}_2$
if and only if it is log-concave.
Furthermore, $h$ is a \textit{P\'olya frequency function}
(and we write $h \in \mathrm{PF}_{\infty} $) if $K(x,y) \equiv h(x-y) $ is
totally positive
of all orders $m$; see, for example, Schoenberg \cite{MR0047732}, Karlin \cite{MR0230102}, and
Marshall, Olkin and Arnold \cite{MR2759813}.
Following Karlin \cite{MR0230102}, we say that $h$ is \textit{strictly
$\mathrm{PF}_{\infty}$} if all the determinants
$\operatorname{det}(H_m) $ are strictly positive.\vspace*{-1pt}

%th2.1 #&#
\begin{theorem}
\label{ChernoffDensityLogConcave}
For each $c>0$ the density $f_{Z_c}(x) = (1/2) g_c(x) g_c(-x)$ is $\mathrm{PF}_{2}$;
that is, log-concave.
\end{theorem}

The Fourier transform in (\ref{GroeneboomFT})
implies that $g_c$ has bilateral Laplace transform (with a slight abuse
of notation)\vspace*{-1pt}
%
%e2.1 #&#
\begin{equation}
\hat{g}_c (z ) = \int \mathrm{e}^{z s} g_c (s) \,
\mathrm{d}s = \frac{2^{1/3}
c^{-1/3}}{\operatorname{Ai} ( (2c^2)^{-1/3} z )} \label{GroeneboomBLT}\vspace*{-1pt}
\end{equation}
for all $z$
such that $\operatorname{Re}(z) > -a_1/(2c^2)^{-1/3}$ where $-a_1$ is the
largest zero of $\operatorname{Ai}(z)$
in $(-\infty, 0)$.

To prove Theorem~\ref{ChernoffDensityLogConcave},
we first show that $g_c$ is $\mathrm{PF}_{\infty}$ by application of the
following two results.

%th2.2 #&#
\begin{theorem}[(Schoenberg, 1951)]
\label{SchoenbergTheorem51}
A necessary and sufficient condition for a (density) function
$g(x)$, $-\infty<x<\infty$, to be a $\mathrm{PF}_{\infty}$ (density) function
is that the reciprocal of its bilateral Laplace transform (i.e., Fourier)
be an entire function of the form
%
%e2.2 #&#
\begin{equation}
\psi(s) \equiv\frac{1}{\hat{g} (s)} = C \mathrm{e}^{-\gamma s^2 + \delta s} s^k \prod
_{j=1}^\infty(1+b_j s)
\exp(-b_j s ), \label{SchoenbergCharacterization}
\end{equation}
where
$C>0$, $\gamma\ge0$, $\delta\in\RR$, $k \in\{0, 1, 2 , \ldots\}
$, $b_j \in\RR$,
$\sum_{j=1}^\infty| b_j |^2 < \infty$.
(For the subclass of densities, the if and only if statement holds for
$1/\hat{g}$ of this form with $\psi(0) =C =1$ and $k=0$.)
\end{theorem}

%pr2.1 #&#
\begin{proposition}[(Merkes and Salmassi)]
\label{ProductRepresentationOfAiry}
 Let $\{- a_k\}$ be the zeros of the Airy function
$\operatorname{Ai}$ (so that $a_k > 0$ for each $k$).
The Hadamard representation of $\operatorname{Ai}$ is given by
\[
\operatorname{Ai} (z) = \operatorname{Ai} (0) \mathrm{e}^{-\nu z} \prod_{k=1}^\infty(1
+ z/a_k) \exp(-z/a_k),
\]
where
\begin{eqnarray*}
 \operatorname{Ai}(0) &= &\frac{1}{3^{2/3} \Gamma(2/3)} = \frac{\Gamma
(1/3)}{3^{1/6} 2 \pi} \approx0.35503,
\\
 \operatorname{Ai}' (0) &=& - \frac{1}{3^{1/3} \Gamma(1/3)} = - \frac{3^{1/6}
\Gamma(2/3)}{2\pi}
\approx-0.25882\quad \mbox{and}
\\
 \nu&=& -\operatorname{Ai}' (0)/\operatorname{Ai}(0) = \frac{3^{1/3} \Gamma(2/3)}{\Gamma(1/3)} = \frac{2 \pi}{3^{1/6} \Gamma(1/3)^2}
\approx0.729011 \ldots .
\end{eqnarray*}
\end{proposition}

\begin{figure}

\includegraphics{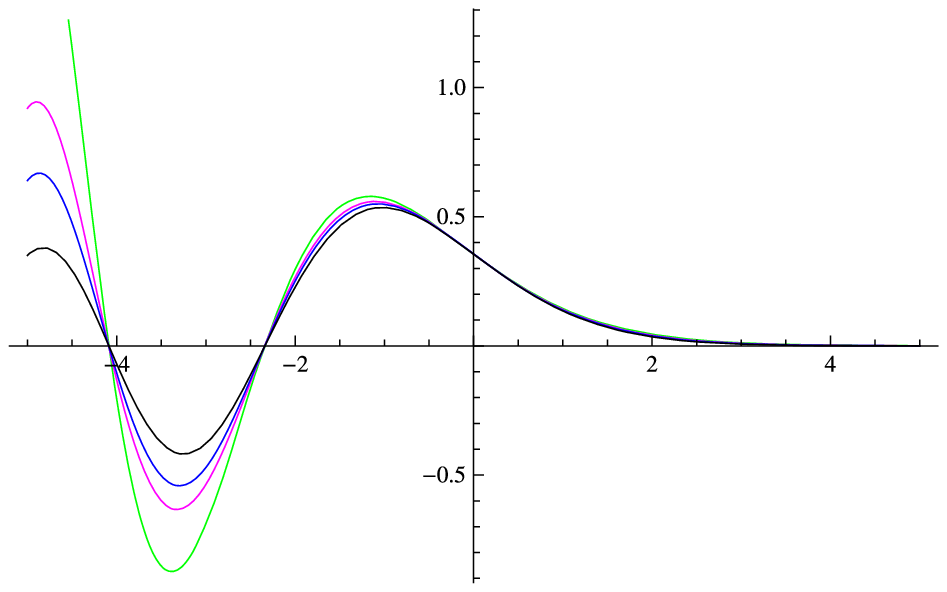}

\caption{Product approximations of $\operatorname{Ai}(x)$.}\label{fig:ApproximateAiry}
\end{figure}

Proposition~\ref{ProductRepresentationOfAiry} is given by Merkes and Salmassi \cite
{MR1624939}; see their Lemma 1, page 211.
This is also Lemma 1 of Salmassi \cite{MR1731663}. Our statement of
Proposition~\ref{ProductRepresentationOfAiry} corrects the constants
$c_1$ and $c_2$ given by
Merkes and Salmassi \cite{MR1624939}.
Figure~\ref{fig:ApproximateAiry}
shows $\operatorname{Ai}(z)$ (black) and $m$ term approximations to $\operatorname{Ai}(z)$ based on
Proposition~\ref{ProductRepresentationOfAiry} with $m=25$ (green),
$125$ (magenta), and $500$ (blue).

%f1 #&#

%pr2.2 #&#
\begin{proposition}
\label{gTotalPositivity}
The functions $t \mapsto g_c (t)$ are in $\mathrm{PF}_{\infty} \subset \mathrm{PF}_2$
for every $c>0$.
Thus, they are log-concave. In fact, $t \mapsto g_c (t)$ is strictly
$\mathrm{PF}_{\infty}$ for every $c>0$.
\end{proposition}

\begin{pf}
By Proposition 2.1,
\begin{eqnarray*}
\operatorname{Ai} \bigl( \bigl(2c^2 \bigr)^{-1/3} z \bigr) & = & \operatorname{Ai} (0)
\mathrm{e}^{-\nu(2c^2)^{-1/3} z} \prod_{j=1}^\infty \biggl(1 +
\frac{z}{(2c^2)^{1/3} a_j} \biggr) \exp \biggl( - \frac{z}{(2c^2)^{1/3} a_j } \biggr)
\\
& = & \operatorname{Ai} (0) \mathrm{e}^{\delta z} \prod_{j=1}^\infty
(1 + b_j z ) \exp ( - b_j z ),
\end{eqnarray*}
which is of the form (\ref{SchoenbergCharacterization}) required in
Schoenberg's theorem with
$k=0$,
%
%e2.3 #&#
%e2.4 #&#
%e2.5 #&#
\begin{eqnarray}
 \delta&= &- \bigl(2c^2 \bigr)^{-1/3} \nu = -
\frac{(3/2)^{1/3}\Gamma(2/3)}{c^{2/3} \Gamma(1/3)} , \label{DeltaSchoenberg}
\\
 C& =& \operatorname{Ai} (0) =1/ \bigl( 3^{2/3} \Gamma(2/3) \bigr) \quad\mbox{and}
\label{LeadConstantSchoenberg}
\\
 b_j &=& \frac{1}{(2c^2)^{1/3} a_j },\qquad j \ge1, \label{B_sSchoenberg}
\end{eqnarray}
where $\{ - a_j \}$ are the zeros of the Airy function $\operatorname{Ai}$. Thus, we
conclude from Schoenberg's theorem
that $g_c$ is $\mathrm{PF}_{\infty}$ for each $c>0$.

The strict $\mathrm{PF}_{\infty}$ property follows from Karlin \cite{MR0230102},
Theorem 6.1(a), page 357: note that in the
notation of Karlin \cite{MR0230102}, $\gamma=0$ and Karlin's $a_i$ is our
$1/a_k$ with
$\sum_k (1/a_k) = \infty$ in view of the fact that
$a_k \sim((3/8)\pi(4k-1))^{2/3}$ via 9.9.6 and 9.9.18, page 18,
Olver \textit{et~al.}~\cite{OlverEtAl:10}.
\end{pf}

Now we are in position to prove Theorem~\ref{ChernoffDensityLogConcave}.

\begin{pf*}{Proof of Theorem~\ref{ChernoffDensityLogConcave}}
This follows from Proposition~\ref{gTotalPositivity}: note that
\[
-\log f_{Z_c} (x) = -\log g_c (x) -\log g_c
(-x) ,
\]
so
\[
w(x) \equiv(-\log f_{Z_c})'' (x) = (-\log
g_c )'' (x) + (-\log g_c
)'' (-x) \equiv v(x) + v(-x) \ge0
\]
since $g_c \in \mathrm{PF}_{\infty} \subset \mathrm{PF}_2$.
\end{pf*}

\textit{Some scaling relations}:
From the Fourier tranform of $g_c$ given above, it follows that
\begin{eqnarray*}
g_c (x) & = & \frac{(2/c)^{1/3}}{2 \pi} \int_{-\infty}^\infty
\frac
{\mathrm{e}^{-\mathrm{i}ux}}{\operatorname{Ai} (\mathrm{i} (2c^2)^{-1/3} u)} \,\mathrm{d}u
\\
& = & \frac{(2/c)^{1/3} (2c^2)^{1/3}}{2 \pi} \int_{-\infty}^\infty
\frac{\mathrm{e}^{-\mathrm{i}v (2c^2)^{1/3} x}}{\operatorname{Ai} (\mathrm{i}v)} \,\mathrm{d}v
\\
& \equiv& 2^{1/6} c^{1/3} g_{2^{-1/2}} \bigl(
\bigl(2c^2 \bigr)^{1/3} x \bigr) .
\end{eqnarray*}
Thus it follows that\vspace*{-1pt}
\[
( \log g_c )'' (x) =
\bigl(2c^2 \bigr)^{2/3} \cdot ( \log g_{2^{-1/2}}
)'' \bigl( \bigl(2c^2 \bigr)^{1/3}
x \bigr) ,
\]
and, in particular,\vspace*{-1pt}
\[
\bigl( \log g_c (x) \bigr)''
|_{x=0} = \bigl(2c^2 \bigr)^{2/3} \cdot ( \log
g_{2^{-1/2}} )'' (x) |_{x=0} .
\]
When $c=1$, the conversion factor is $2^{2/3}$. Furthermore we compute
\begin{eqnarray*}
f_{Z_c} (t) & = & \tfrac{1}{2} g_c (t)
g_c (-t) = \tfrac{1}{2} 2^{1/3} c^{2/3}
g_{2^{-1/2}} \bigl( \bigl(2c^2 \bigr)^{1/3}t \bigr)
g_{2^{-1/2}} \bigl(- \bigl(2c^2 \bigr)^{1/3} t \bigr)
\\
& \equiv& c^{2/3} f_1 \bigl( c^{2/3} t \bigr),
\end{eqnarray*}
where\vspace*{-1pt}
\begin{eqnarray*}
f_1 (t) & \equiv & f_{Z_1} (t) = \tfrac{1}{2}
g_1 (t) g_1 (-t)
\\[-1pt]
& = & \tfrac{1}{2} 2^{1/3} g_{2^{-1/2}} \bigl(2^{1/3}
t \bigr) g_{2^{-1/2}} \bigl(- 2^{1/3} t \bigr) .
\end{eqnarray*}
Thus we see that\vspace*{-1pt}
\[
Z_c \stackrel{d} {=} c^{-2/3} Z_1
\]
for all $c>0$.

\begin{figure}[t]

\includegraphics{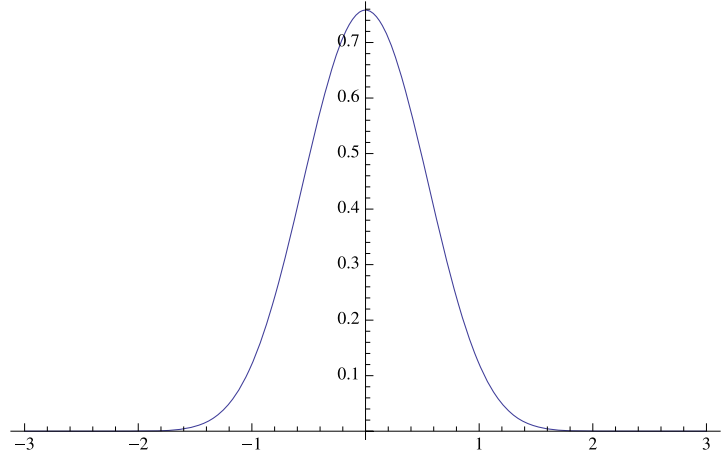}

\caption{The density $f_Z$.}\label{fig:ChernoffDensity}
\end{figure}

\begin{figure}[b]

\includegraphics{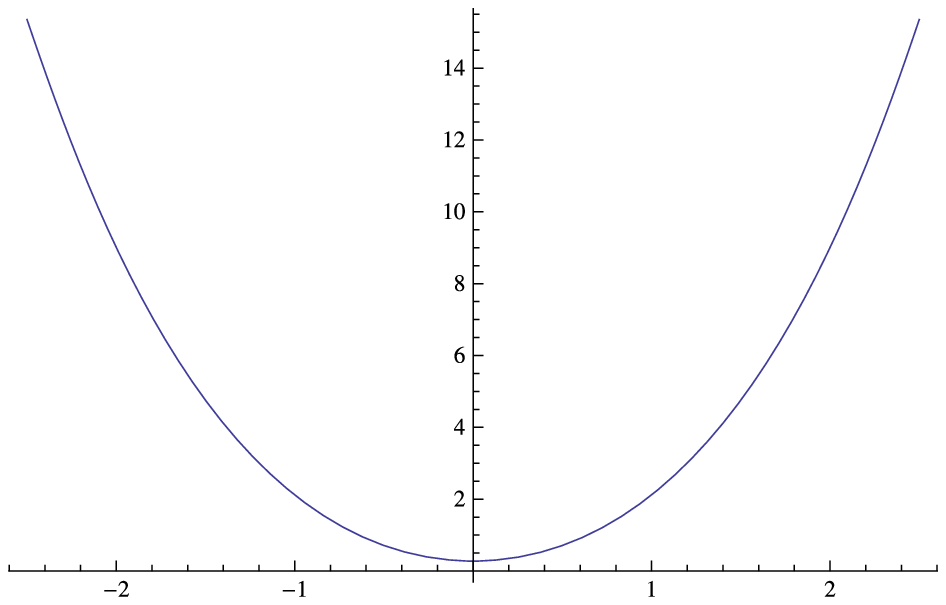}

\caption{$-\log f_Z$.}\label{fig:MinusLogChernoff}
\end{figure}

\begin{figure}

\includegraphics{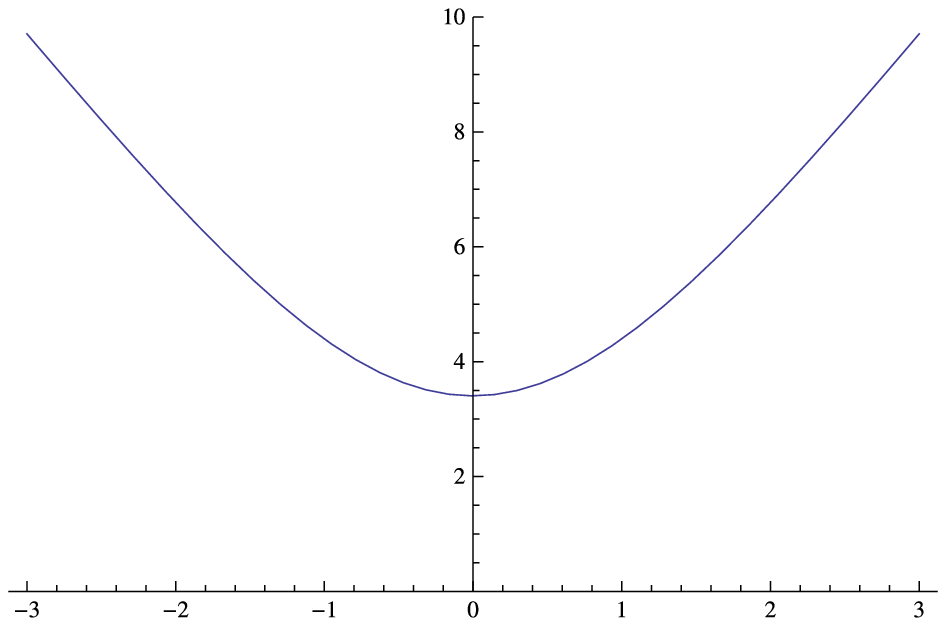}

\caption{$(-\log f_Z)''$.}\label{fig:SecondDerivativeMinusLogChernoff}
\end{figure}

Figure~\ref{fig:ChernoffDensity} gives a plot of $f_Z$;
Figure~\ref{fig:MinusLogChernoff} gives a plot of $-\log f_Z$; and
Figure~\ref{fig:SecondDerivativeMinusLogChernoff} gives a plot of
$(-\log f_Z)'' $.

%f2 #&#

%f3 #&#

%f4 #&#

If we use the inverse Fourier transform to
represent $g$ via (\ref{GroeneboomFT}), and then calculate directly,
some interesting correlation type inequalities involving the Airy
kernel emerge.
Here is one of them.

Let $h(u) \equiv1/| \operatorname{Ai} (\mathrm{i}u)| \sim2 \sqrt{\pi} u^{1/4} \exp(-
(\sqrt{2}/3) u^{3/2})$ as $u \rightarrow\infty$
by Groeneboom \cite{MR981568}, page~95. We also define
$\varphi(u,x/2) = \operatorname{Re} (\mathrm{e}^{\mathrm{i}ux/2} \operatorname{Ai} (\mathrm{i}u)) h(u)$ and $\psi(u,x/2) = \operatorname{Im}
( \mathrm{e}^{\mathrm{i}ux/2}\times \operatorname{Ai} (\mathrm{i}u)) h(u)$.

%co2.1 #&#
\begin{corollary} With the above notation,
\begin{eqnarray*}
&&\int_0^{\infty} \sin^2 (uy)
\varphi(u,x) h(u) \,\mathrm{d}u \cdot \int_0^\infty
\cos^2 (uy) \varphi(u,x) h(u) \,\mathrm{d}u
\\[-1pt]
&&\quad{} + \int_0^\infty\sin(uy) \cos(uy) \psi(u,x)
h(u) \,\mathrm{d}u \ge0 \qquad\mbox{for all } x, y \in\RR.
\end{eqnarray*}
\end{corollary}

%s3 #&#
\section{Is Chernoff's density strongly log-concave?}
\label{sec:StrongLogConcave}

From Rockafellar and Wets \cite{MR1491362} page 565,
$h \dvt \RR^d \rightarrow\overline{\RR}$ is \textit{strongly convex} if
there exists a constant $c>0$ such
that
\[
h \bigl(\theta x + (1-\theta)y \bigr) \le\theta h(x) + (1-\theta) h(y) -
\tfrac
{1}{2} c \theta(1-\theta) \| x - y \|^2
\]
for all $x,y \in\RR^d$, $\theta\in(0,1)$. It is not hard to show
that this is equivalent to convexity of
\[
h(x) - \tfrac{1}{2} c \| x \|^2
\]
for some $c>0$. This leads (by replacing $h$ by $-\log f$)
to the following definition of \textit{strong log-concavity} of a
(density) function:
$f \dvtx \RR^d \rightarrow\overline{\RR}$ is strongly log-concave if
and only if
\[
-\log f(x) - \tfrac{1}{2} c \| x \|^2
\]
is convex for some $c >0$.
Defining $-\log g(x) \equiv-\log f(x) -(1/2)c \| x \|^2$, it is easily seen
that~$f$ is strongly log-concave if and only if
\[
f(x) = g(x) \exp \bigl(-(1/2) c \| x \|^2 \bigr)
\]
for some $c >0$ and log-concave function $g$. Thus if $f \in C^2 (\RR^d)$, a sufficient condition
for strong log-concavity is: $\operatorname{Hess} (-\log f) (x) \ge c I_d$ for
all $x \in\RR^d$ and some $c>0$
where $I_d $ is the $d\times d$ identity matrix.

Figure~\ref{fig:SecondDerivativeMinusLogChernoff} provides compelling evidence
for the following conjecture concerning strong log-concavity of
Chernoff's density.

%co3.1 #&#
\begin{conjecture}
\label{ChernoffDensityLogConcaveStrongForm}
Let $Z_1$ again be a ``standard'' Chernoff random variable.
Then for $\sigma\ge\sigma_0 \approx 0.541912 \ldots = (-(\log f_{Z_1}
(z))'' |_{z=0})^{-1/2} $
the density $f_{Z_1}$ can be written as
\[
f_{Z_1} (x) = \rho(x) \frac{1}{\sigma} \varphi(x/\sigma),
\]
where $\varphi(x) = (2\pi)^{-1/2} \exp(-x^2/2)$ is the standard
normal density
and $\rho$ is log-concave.
Equivalently, if $c \ge\sigma_0^{3/2} \approx0.398927\ldots,$ then
\[
f_{Z_c} (x) = \tilde{\rho} (x) \varphi(x),
\]
where $\tilde{\rho}$ is log-concave.
\end{conjecture}

\begin{pf} (Partial) Let $w \equiv(-\log f_{Z_c})''$ and $v \equiv
(-\log g_c)''$.
Then
\[
w(t) = v(t) + v(-t) \ge 2 v(0) = w(0) > 0
\]
is implied by convexity of $v$ and strict positivity of $w(0)$.
Thus, we want to show that $v^{(2)} = (-\log g_c)^{(4)} \ge0$.

To prove this, we investigate the normalized version of $g_c$ given by
$\widetilde{g}_c (x) = g_c (x) \operatorname{Ai}(0) /  (2/c)^{1/3} = g_c (x)/ \int g_c
(y)\,\mathrm{d}y$ so
that $\int\widetilde{g}_c (x) \,\mathrm{d}x =1$.
Suppose that $b_i $ is given in (\ref{B_sSchoenberg}), and
let $X_i \sim\operatorname{Exp}(1/b_i)$ be independent exponential random
variables for $i=1,2,\ldots.$
Since $\sum_{i=1}^\infty b_i^2 < \infty$, the random variable $Y_0 =
\sum_{i=1}^\infty(X_i - b_i)$ is finite almost surely
(see, e.g., Shorack \cite{MR1762415}, Theorem 9.2, page 241)
and the Laplace transform of $-(\delta+Y_0)$ is given by
\begin{eqnarray*}
\varphi(s) & \equiv& \mathrm{e}^{-\delta s} E \mathrm{e}^{-s Y_0} = \exp( - \delta s)
\cdot\frac{1}{ \prod_{i=1}^\infty(1 + b_i s)
\mathrm{e}^{-b_i s} }
\\
& = & \frac{1}{\mathrm{e}^{ \delta s} \cdot \prod_{i=1}^\infty(1 + b_i s)
\mathrm{e}^{-b_i s} } ,
\end{eqnarray*}
exactly the form of the Laplace transform of $g$ implicit in the proof of
Proposition~\ref{gTotalPositivity}, but without
the Gaussian term. Thus, we conclude that $\tilde{g}_c$ is the
density of $Y \equiv - \delta- Y_0 = -\delta- \sum_{j=1}^\infty
(X_j - b_j )$.

Now let $\lambda_i = 1/b_i$ for $i \ge1$. Thus, $X_i \sim\operatorname
{Exp}(\lambda_i)$.
A closed form expression for the density of $Y_m \equiv\sum_{i=1}^m
X_i$ has been given by
Harrison \cite{MR1039185}. %Harrison (1990).
From Harrison's theorem 1, $Y_m$ has density
%
%e3.1 #&#
\begin{equation}
f_m (t) = \sum_{j=1}^m
\lambda_j \exp( - \lambda_j t) \prod
_{i
\not= j} \frac{\lambda_i}{\lambda_i - \lambda_j } . \label{HarrisonDensityFormula}
\end{equation}
If we could show that $v_m (t) \equiv(-\log f_m)'' (t)$ is convex,
then we would be done!
Direct calculation shows that this holds for $m=2$, but our attempts at
a proof for general $m$
have not (yet) been successful.
On the other hand, we know that for $t\ge0$,
\[
w(t) = v(t) + v(-t) \ge v(t) \ge v(0)>0
\]
if $v$ satisfies $v(t) \ge v(0)$ for all $t\ge0$, so we would have
strong log-concavity with the constant~$v(0)$.
\end{pf}

%s4 #&#
\section{Discussion and open problems}
\label{sec:OpenProbs}

Log-concavity of Chernoff's density implies that the peakedness results of
Proschan \cite{MR0187269} and Olkin and Tong \cite{MR927147} apply. See also Marshall and Olkin \cite
{MR552278}, page 373,
and Marshall, Olkin and Arnold~\cite{MR2759813}.

Note that the conclusion of Conjecture~\ref
{ChernoffDensityLogConcaveStrongForm} is
exactly the form of the hypothesis of the inequality of Harg{\'e} \cite
{MR2095937} and of
Theorem 11, page 559, of Caffarelli \cite{MR1800860};
see also Barthe \cite{MR2275657}, Theorem~2.4, page 1532.
Wellner \cite{Wellner:12} shows that the class of strongly log-concave
densities is closed under convolution, so in particular if
Conjecture~\ref{ChernoffDensityLogConcaveStrongForm}
holds, then the sum of two independent Chernoff random variables is
again strongly
log-concave.

Another implication is that a theorem of Caffarelli \cite{MR1800860} applies: the
transportation map $T = \nabla\varphi$ is a contraction.
In our particular one-dimensional special
case, the transportation map $T $ satisfying $T(X) \stackrel{d}{=} Z$
for $X \sim N(0, 1)$ is just the solution of $\Phi(z) = F_Z (T(z))$, or
equivalently $T(z) = F_{Z}^{-1} ( \Phi(z))$. This function is
apparently connected
to another question concerning convex ordering of $F_Z$ and $\Phi
(\cdot)$ in the sense of
van Zwet \cite{MR0176511}; see also van Zwet \cite{MR0175217}: is $T^{-1} (w) =  \Phi^{-1} ( F_Z (w))$ convex
for $w >0$?

As we have seen above,
Chernoff's density has the symmetric product form (\ref
{SymmetricProductFormChernoff})
where $g$ has Fourier transform given in (\ref{GroeneboomFT}).
In this case, we know from Section~\ref{sec:LogConcave} that $g \in
\mathrm{PF}_{\infty}$.

As is shown in the longer technical report version of this paper Balabdaoui and
  Wellner \cite
{Balabda-Wellner:12},
it follows from the results of Bondesson \cite{MR1224674,MR1481175} that
the standard normal density $\phi$ can be written
in the same structural form as that of Chernoff's density (\ref
{SymmetricProductFormChernoff});
that is:
%
%e4.1 #&#
\begin{equation}
\phi(z)  = \tfrac{1}{2} g(z) g(-z) ,\label{SecondIdentity}
\end{equation}
where now
%
%e4.2 #&#
\begin{eqnarray}\label{SecondFormGforStdGauss}
g(z) & \equiv& (2 / \pi)^{1/4} \exp(z) \exp \biggl( \int
_0^\infty \log \biggl( \frac{\mathrm{e}^s + 1}{\mathrm{e}^s +\mathrm{e}^z} \biggr)
\,\mathrm{d}s \biggr)
\nonumber
\\[-8pt]
\\[-8pt]
\nonumber
& = & (2 / \pi)^{1/4} \exp \biggl(\frac{\pi^2}{12}+z - \int
_0^{\mathrm{e}^z} \frac{\log(1+t)}{t} \,\mathrm{d}t \biggr)
\end{eqnarray}
is log-concave, integrable, and $g \in\log(\mathrm{HM}_{\infty})$,
the log-transform (in terms of random variables) of the Hyperbolically
Completely Monotone
class of Bondesson \cite{MR1224674,MR1481175}. %Bondesson (xx,yy)
Two natural questions are:
(a) Does the function $g$ in (\ref{SymmetricProductFormChernoff})
satisfy $g \in\log(\mathrm{HM}_{\infty})$?
(b) Does the function $g$ in (\ref{SecondFormGforStdGauss}) satisfy $g
\in \mathrm{PF}_{\infty}$?

A further question remaining from Section~\ref{sec:StrongLogConcave}:
Is Chernoff's density strongly log-concave?

A whole class of further problems involves replacing the (ordered)
convex cone
$K_n$ in Section~\ref{sec:TwoLimitTheorems} by
the convex cone $\widetilde{\cal K}_n$ corresponding to a convexity
restriction as in Section 2 of
Groeneboom, Jongbloed and
  Wellner \cite{MR1891742}. In this latter case, the limiting distribution
depends on an ``invelope'' of
the integral of a two-sided Brownian motion plus a polynomial drift as
follows: it is the density of
the second derivative at zero of the ``invelope''. See Groeneboom, Jongbloed and
Wellner \cite{MR1891742,MR1891741} for further
details and Balabdaoui, Rufibach and
  Wellner \cite{MR2509075} for another convexity related shape
constraint where this limiting distribution occurs.
However, virtually nothing is known concerning the analytical
properties of this distribution.

\section*{Acknowledgements} We owe thanks to Guenther Walther for pointing
us to Karlin \cite{MR0230102} and
Schoenberg's theorem. We also thank Tilmann Gneiting for several
helpful discussions.
Supported in part by NSF Grants DMS-08-04587 and
DMS-11-04832, by NI-AID Grant 2R01 AI291968-04, and by the Alexander
von Humboldt Foundation.

% imsref loaded by akundreckaite, 2013-01-11 15:15:32
% imsref loaded by akundreckaite, 2013-01-14 09:31:21
% imsref loaded by akundreckaite, 2013-01-14 09:43:04

%suskaldyti doi

\printhistory

\end{document}